\documentclass[12pt,thmsa]{article}
\usepackage{amsfonts}
\usepackage{amssymb}
\newtheorem{teo}{Theorem}[section]

\newtheorem{obs2}[teo]{Remark}

\newtheorem{tea}{Theorem}[subsection]

\newtheorem{no2}[teo]{Note}

\newtheorem{no3}[tea]{Note}

\newcommand{\F}{{\mathbb{F} }}
\newcommand{\Q}{{\mathbb{Q} }}

\title{Appendix to: The level 1 weight 2 case of Serre's conjecture - a strategy for a proof}
\author{Luis V. Dieulefait\\ Universitat de Barcelona
  \thanks{
e-mail: ldieulefait@ub.edu
}}

\begin{document}

\maketitle

\section{A fact, a question, and its answer}
(we follow the notation of the preprint)

 As explained in the preprint, to complete the proof our task is to show that, among the minimal (modular) Barsotti-Tate p-adic deformations of 
 $\hat{\rho}|_{G_F}$ there is at least one that can be extended to $G_\Q$.\\

The field $F$ is known to be unramified at $p$ and we can choose it (as explained in the preprint, using solvable base-change and Sylow theorems) so that its degree over $\Q$ is prime to $p$.\\

 With these conditions, as we already explained, there is a one-to-one correspondence between minimally ramified Barsotti-Tate deformations of $\hat{\rho}$ and  those minimally ramified Barsotti-Tate deformations of $\hat{\rho}|_{G_F}$ that can be extended to the full $G_\Q$. (CC)\\
 
 Let us call $R$ the universal deformation ring of minimally ramified (Barsotti-Tate at $p$) deformations of $\hat{\rho}$, and let $R'$ denote a similar minimal universal deformation ring, but of $\hat{\rho}|_{G_F}$.\\
 We know from the results of Taylor that $R'$ is a complete intersection ring.\\
  We will also need the following fact (proved by Bockle and Ramakrishna):\\
 
 FACT  (FF): Let $W$ be the Witt ring of $\F_q$, then  $R$ is an $W$-algebra of the type:
 
 $$  W [[ X_1, .... X_r]] / (f_1, ...., f_s) $$
 
 with $r \geq s$.\\
 
 Now, let us answer the question: Why in the complete intersection ring $R'$ there is at least one p-adic deformation that survives when we descend to $G_\Q$?\\
 Well, if none of the minimal deformations of $\hat{\rho}|_{G_F}$  descends to $G_\Q$, then using the correspondence (CC) we see that $R$ is too small to match with (FF).\\
 This concludes the proof, so existence of minimal deformations, and Serre's conjecture in the level 1 weight 2 case follow.\\

\end{document}